\documentclass[12pt]{amsart}

\usepackage{amssymb,amsmath,amsthm,epsfig,amscd,amssymb,latexsym,graphicx, epic, eepic}
\usepackage{color}
\usepackage{hyperref}

\numberwithin{equation}{section}

\def\H{\mathcal H}

\newtheorem{definition}{Definition}[section]
\newtheorem{theorem}{Theorem}[section]

\newtheorem{corollary}[theorem]{Corollary}

\newcommand{\R}{{\mathbb R}}

\newcommand{\C}{{\mathbb C}}

\newcommand{\N}{{\mathbb N}}

\def\hat{\widehat}

\def\bfo{\begin{eqnarray*} }
\def\efo{\end{eqnarray*} }
\def\ba{\begin{eqnarray*} }
\def\ea{\end{eqnarray*} }
\def\beq{\begin{eqnarray}}
\def\eeq{\end{eqnarray}}
\def\supp{\hbox{supp}\,}

\def\p{\partial}

\def\p{\partial}

\def\R{\mathbb R}
\def\oo{\it o}

\begin{document}

\title{Inverse Problems for Screens}

\author{Emilia Bl{\aa}sten, Petri Ola and Lassi P\"aiv\"arinta}



\begin{abstract}

\smallskip We study the inverse scattering from a screen with using only one incoming time--harmonic plane wave but with measurements of the scattered wave done at all directions. Especially we focus on the \(2D\)--case i.e. (inverse) scattering from an open bounded smooth curve. Besides the inverse scattering problem we also study the inverse electrostatic problem. We then show that one Cauchy--data of any continuous and bounded function vanishing on the screen and harmonic outside it, determines the screen uniquely.

\end{abstract}

\maketitle

\section{Direct Scattering Problem} 
\smallskip

An open set  \(\Omega \subset \R^n\), \(n\geq 2\), is called a {\em Dirichlet--obstacle} if \(\Omega\) is bounded and the complement \(\R^n \setminus \overline \Omega\) is connected. The corresponding direct scattering problem is then the following: For a fixed wavenumber \(k^2>0\) let the incident wave \(u_i\) be a solution of the Helmholtz--equation  \((\Delta + k^2) u_i = 0\) in \(\R^n\). Find the {\em total field} \(u\in C^2(\R^n \setminus \overline \Omega) \cap C(\R^n \setminus \Omega)\) such that 
\begin{itemize}
\item\((\Delta + k^2) u = 0 \quad \text{in \(\R^n \setminus \overline \Omega\)}\)

\item The {\em scattered field} \(u_s = u - u_i\) is outgoing in the sense that it satisfies Sommerfeld's Radiation condition
\[
(\partial_ r u_s - ik u_s)(x)  = \oo (|x|^{-(n-1)/2}) \quad \text{as \(|x| \to \infty\) uniformly in \(x/|x|\in S^{n-1}\)}
\]
\item
The total field vanishes on \(\partial \Omega\), i.e. the scattered field satisfies \(u_s = - u_i\) on \(\partial \Omega\).
\end{itemize}
\smallskip

\noindent This is a classic problem and is understood very well, even if \(\Omega\) has minimal regularity (see  for example \cite{McL}, \cite{NPT}). The situation is less well understood if \(\Omega\) is replaced by a compact surface \(S\), i.e. a {\em screen} as defined below:

\begin{definition}
An \((n-1)\)--dimensional screen is a smooth compact orientable sub\-mani\-fold \(S\) of \(\R^{n}\) with nonempty boundary \(\partial S\).
\end{definition}

\noindent Hence, a two--dimensional screen is a smooth compact hypersurface of \(\R^3\), and a one--dimensional screen is a compact smooth curve. Both with nonempty boundaries of course. The scattering problem for \((n-1)\)--dimensional screens that we are going to study is then  as follows: For \(u_i\) as above, find the scattered field \(u_s\) such that the total field \(u = u_s + u_i\) solves
\begin{enumerate}
\item \((\Delta + k^2) u = 0 \quad \text{in \(\R^n \setminus S\)}\), \label{screen scattering 1}

\item The {\em scattered field} \(u_s = u - u_i\) is outgoing i.e. it satisfies Sommerfeld's Radiation condition \label{screen scattering 2}
\[
(\partial_ r u_s - ik u_s)(x)  = \oo (|x|^{-(n-1)/2}) \quad \text{as \(|x| \to \infty\) uniformly in \(x/|x|\in S^{n-1}\)}
\]

\item The total field vanishes on \(S\), i.e. the scattered field satisfies \(u_s = - u_i\) on \(S\). \label{screen scattering 3}
\end{enumerate}
\smallskip

\noindent Physically this describes acoustic time--harmonic scattering from a sound--soft screen. For the solvability in dimension three and in the case when \(S\) can be smoothly embedded to a boundary of a smooth bounded domain, see for example \cite{St}.
\medskip

\section{Inverse Scattering} If \(u_s\) is any outgoing solution of the Helmholtz--operator \(\Delta + k^2\) then at infinity it has asymptotics
\[
u_s(x) = \frac{e^{ik|x|}}{|x|^{(n-1)/2}}u_\infty (\hat x, k) + \oo (|x|^{-(n-1)/2}), \quad \hat x = x/|x|.
\]
The uniquely determined function \(u_\infty\) is the far--field pattern of \(u_s\). In the case when \(u_s\) is scattered wave corresponding to the incoming plane wave \(u_i (x) = e^{ik\langle \theta, x\rangle}\) we denote it by \(u_\infty (\hat x, \theta, k)\). Let's give the following rather general formulation, this will be narrowed down considerably in a moment. \smallskip

\noindent {\bf Inverse Scattering Problem (ISP)}: For a fixed \(k>0\) determine the obstacle \(\Omega\) (or the screen \(S\)) from the knowledge of \(u_\infty (\hat x, \theta, k)\) on a given set of values of \((\hat x, \theta) \in S^{n-1} \times S^{n-1}\). 
\smallskip

A classical result of Schiffer (see for example \cite{CK-II,LP,T2}) says that given any countably infinite set \(\theta_j\), \(j\in \N\), of incoming directions,  the values \(u_\infty (\hat x, \theta_j , k)\) for all measurement directions \(\hat x \in S^{n-1}\) determine \(\Omega\) uniquely. This fact holds also for sound--soft screens \cite{K-fixed-freq}. {\em The Schiffer's problem} then asks whether this holds with just one fixed incoming direction \(\theta\), i.e whether the measurement
\[
M_\theta = \{u_\infty(\hat x, \theta, k); \, \hat x \in S^{n-1}\}
\] 
is enough to determine the obstacle \(\Omega\). 

\smallskip
We propose the following: 
\medskip

\noindent{\em Conjecture:} A single measurement set \(M_\theta = \{u_\infty(\hat x, \theta, k); \, \hat x \in S^{n-1}\}\) determines a sound--soft screen uniquely.
\medskip

For previous results on inverse scattering from screens, see \cite{AH,K-analytic}.

\medskip

\section{Screens in Two Dimensions}
 In \(\R^2\) screens are arcs, and hence geometrically and analytically more simple than in higher dimensions. So, from now on we will consider \eqref{screen scattering 1}--\eqref{screen scattering 3} with \(S\) replaced by a simple differentiable arc \(\Gamma\). An integration of parts immediately  gives that the solution of the scattering problem - if it exists - must be 
 \begin{equation}\label{equ}
   u_s (x) = - \frac{i}{4} \int_\Gamma H_0^{(1)} (k|x-y|) \rho (y) \, ds(y),
 \end{equation}
where the density \(\rho\) is equal to the jump \([\partial_\nu u_s]\) of the normal derivative across \(\Gamma\) and \(H_0^{(1)}\) is the Hankel--function of the first kind and order zero. Note that  \([\partial_\nu u_s]\) is independent of the direction of the normal \(\nu\) of \(\Gamma\) chosen. Taking limits on \(\Gamma\) then gives that \(\rho\) must solve an integral equation
\begin{equation}
  - u_i (x) = -\frac{i}{4} \int_\Gamma H_0^{(1)} (k|x-y|) \rho (y) \, ds(y), \quad x\in \Gamma.
\end{equation}
This is a convolution equation\footnote{In fact, the right hand side defines a pseudo--differential operator of order -1 on \(\Gamma\) when \(\Gamma\) is \(C^\infty\)--smooth. However, the presence of boundary creates additional problems when analysing solvability.} on \(\Gamma\). 
\bigskip
\smallskip

Instead of the scattering problem above it will be easier to study the corresponding electrostatic question, i.e. when \(k = 0\). Hence we wish to determine the screen \(\Gamma\) from electrostatic measurements. In this case the inverse--problem can be formulated as whether one can uniquely determine a (super) conductive crack in a homogenous body from electrostatic measurements. More precisely, assume that \(\Gamma\) is a screen in \(\R^2\), and for some constant \(C\) we have
\begin{eqnarray}\label{static screen}
\Delta u & = & 0 \quad \text{in \(\R^2 \setminus \Gamma\)}, \\
u|_\Gamma & = & 0\\
u & = & u_0 + C \quad \text{in \(\R^2 \setminus \Gamma\)},\label{static screen end}
\end{eqnarray}
where \(u_0\) has finite energy in the sense that
\[
\int_{\R^2 \setminus \Gamma} |\nabla u_0|^2 \, dx + \int_\Gamma |u_0|^2 \, dx <\infty.
\]
We consider the following:
\smallskip

\noindent{\em Electrostatic Schiffer's Problem}: Let \(B\) be a closed ball in \(\R^2\) containing \(\Gamma\) in its interior. Given a single solution \(u\) of \eqref{static screen}--\eqref{static screen end}, determine \(\Gamma\) from the Cauchy--data \((u|_{\partial B}, \partial_r u|_{\partial B})\) of \(u\) on \(\partial B\). 
\smallskip

\noindent The essential restriction here is that we have only one single solution at our disposal. This is in marked contrast to the version of classical Calder\'on's problem for screens, where the data would be the Cauchy--data on \(\partial B\) of all functions harmonic in the complement of \(\Gamma\). In \cite{FV} the authors studied a similar problem but with data consisting of Cauchy--data of a pair of solutions in a bounded domain containing an insulating screen in its interior. In that case at least two measurements are needed. In \cite{K-integral,AK} the authors proved uniqueness and numerical recovery for determining an unknown Dirichlet obstacle from one electrostatic measurement on the known boundary. In \cite{AK} they briefly mention that their method could work in principle for open arcs (screens) but we could not find further evidence about this in the literature. Nevertheless, the study of screens is interesting because of the sinularities generated by their boundary.

\section{Unique determination for the Electrostatic Problem} Let \(u\) be a solution of \eqref{static screen}--\eqref{static screen end} and denote
\[
D_\Gamma (u)  = (u|_{\partial B}, \partial_r u|_{\partial B}).
\]

Our main technical result is the following:

\begin{theorem}\label{sing}
If \(u\) is a non--vanishing solution of \eqref{static screen}--\eqref{static screen end}, then \(u\) is singular at boundary points  of \(\Gamma\). 

\end{theorem}

\noindent More precisely we can actually prove that near the endpoints \(b_\pm\) of \(\Gamma\) the solution of \(u\) behaves like
\[
u(x) \sim A_\pm {\rm dist}\, (x, b_\pm)^{-1/2},
\]
where \(A_\pm \not = 0\). As a corollary of this we get the following unique determination result:

\begin{corollary}\label{unique} (Bl{\aa}sten, Ola, P\"aiv\"arinta (2024), in preparation) Let \(\Gamma_1\) and \(\Gamma _2\) be two screens contained in the interior of a closed ball \(B\subset \R^2\) and assume \(u_1\) and \(u_2\) are a pair of corresponding non--vanishing solutions of \eqref{static screen}--\eqref{static screen end}. If \(D_{\Gamma_1} (u_1) = D_{\Gamma_2} (u_2)\) then \(\Gamma_1 = \Gamma_2\).
\end{corollary}

Next we aim to give outlines of the proofs of these two results, and start by showing how Theorem \ref{sing}  implies Corollary \ref{unique}. Assume the contrary, i.e. that \(\Gamma_1 \not  = \Gamma_2\). The endpoints of both of these curves must be in the closure of the connected component of \(B\setminus (\Gamma_1 \cup \Gamma_2)\) containing \(\partial B\). Since \(D_{\Gamma_1} (u_1) = D_{\Gamma_2} (u_2)\) Holmgren's uniqueness Theorem implies that \(u_1 = u_2\) there. Especially, if \(b\) is an endpoint of say \(\Gamma_1\) that does not lie on \(\Gamma_2\), then \(u_2\) must be real--analytic in its neighbourhood, but this contradicts the fact that \(u_1\) is singularr at \(b\). Thus the endpoints of \(\Gamma_1\) and \(\Gamma_2\) must coincide. Consider next the closed curve\footnote{We need to define suitable parametrisations for both of these curves} \(\gamma = \Gamma_1 \Gamma_2^{-1}\) and assume that \(U\) is a connected non--empty  interior component of \(B\setminus \gamma\). The boundary of \(U\) consists of components of \(\Gamma _1\) and \(\Gamma_2\) on which \(u_1 = u_2\) and since \(u_i\) vanishes on \(\Gamma_i\), both \(u_1\) and \(u_2\) are harmonic in \(U\) with vanishing boundary values, hence by maximum principle they must be identically zero outside \(\Gamma _1\) and \(\Gamma_2\) respectively. This is a contradiction so \(\Gamma_1 = \Gamma_2\).
\medskip

\section{Outline of the Proof the Theorem \ref{sing}}
Consider first the case when \(\Gamma\) is the straight line segment from \(-1\) to \(1\), i.e \(
\Gamma = [-1,1]\).
In the electrostatic case the Hankel function in \eqref{equ} is replaced by the fundamental solution of the Laplacian, namely the logarithm. Then integration by parts gives
\[
2\pi u(z) = \int_{-1}^1 \ln |z-s| \, \rho (s) \, ds, \quad z\in \C\setminus \Gamma,
\]
and \(\rho\) satisfies
\begin{equation}\label{curve equ}
 \int_{-1}^1 \ln |t-s| \, \rho (s) \, ds = {\rm const} \neq 0, \quad t\in \Gamma.
\end{equation}
We have to analyse the behaviour of \(u(z)\) when \(z\to \pm 1\). Physically \(u\) represents the electric potential outside \(\Gamma\) and \(\rho\) is the charge density on \(\Gamma\). Differentiating \eqref{curve equ} gives
\[
{\rm PV}\, \int_{-1}^1 \frac{1}{t-s} \rho(s) \, ds = 0, \quad t\in [-1,1].
\]
i.e \(\rho\) belongs to the kernel of the {\em local Hilbert--transform} \(\H_I\) on \(\Gamma\). In the case when \(I\) is replaced by the real line, \(\H = \H_{\R}\) is an isomorphism on \(L^p(\R)\), \(1<p<\infty\), and \(\H^2 = -I\). Also, if \(f = u + iv\) with \(u, \, v\in L^p(\R)\) real valued, a classical Theorem of Riesz says that function \(f\) has an analytic continuation to upper half--space \(\C_+\) if and only if \(v = \H (u)\). Also, if this holds, the analytic extension of \(f\) will belong to the Hardy--space \(H^p(\C_+)\).
Recall that \(F:\C_+ \to \C\) belongs to \(H^p(\C_+)\) if the translates \(F(\cdot + iy)\) belong to \(L^p(\R)\) for all \(y>0\) with a uniform bound their \(L^p\)--norms. Also, every element of \(H^p(\C_+)\) has a non--tangential limit on \(\R\) and this limit belongs to \(L^p(\R)\). However, we have to be careful about the supports when considering \(\H\) instead of \(\H_I\), since \(\H_I^2 \neq -I\).
\smallskip

With this more general point of view we can formulate our question as follows: Assume that \(\rho\not = 0\) satisfies 
\begin{equation}\label{supp 1}
\supp (\rho) \subset I,
\end{equation}
\begin{equation}\label{supp 2}
\supp (\H (\rho)) \subset \R \setminus I.
\end{equation}
Is \(\p_x u = \H (\rho)\) then singular at endpoints of \(I\)? It turns out that it always is, and to see this using the analytic extension is extremly useful. Namely consider the function 
\[
f(z) = (1- z^2)^{-1/2}, \quad z\in \C\setminus i\R_-.
\]
This belongs to \(H^p(\C_+)\) for \(1<p<\infty\) with boundary values
\[
f(t) = \chi_I (t) \, (1- t^2)^{-1/2} + i \chi_{\R \setminus I} (t^2 - 1)^{-1/2}, t\in \R.
\]
i.e if we define \(\rho (t) =  \chi_I (t) \, (1- t^2)^{-1/2}\) we have \(\H(\rho)(t) = \chi_{\R \setminus I} (t^2 - 1)^{-1/2}\), i.e exactly of the form that we are looking for. It turns out all non--zero solutions of satisfying \eqref{supp 1} -- \eqref{supp 2} are of this form and hence they have a non--vanishing inverse square root singularity at the end points of the curve. 
\smallskip

The general case can be reduced to this using the Riemann--mapping theorem on the complement of \(\Gamma\) in the Riemann--sphere \(\hat \C\setminus \Gamma\). Hence there is biholomorphic mapping \(\psi: \hat \C\setminus \Gamma \to \hat \C\setminus [-1,1]\), and using Caratheodory's Prime End Theorem this can be extended to a \(C^1\)--map \(\psi: \Gamma \to [-1, 1]\). Hence singularities on \([-1, 1]\) carry over to \(\Gamma\) by the inverse \(\psi ^{-1}\). 
\medskip

\section{Final Remarks in the Higher Dimensional Case}
Consider a {\em flat screen} \(\Sigma\) in \(\R^3\), i.e we assume that \(S\) is a bounded 2--manifold with boundary and there is a hyperplane \(\pi\) such that \(\Sigma \subset \pi\). In this case the answer to Schiffer's Problem is positive for both in the case of acoustic scattering (\cite{BPS}) and electromagnetic scattering (\cite{OPS}). In both of these cases one can  prove that both the screen \(\Sigma\) and the supporting hyperplane \(\pi\) are uniquely determined by measuring everywhere the non--vanishing scattered wave corresponding to one incoming plane wave. Also, for an exposition on the connection between the Hilbert--, Mellin -- and Fourier--transforms see \cite{BPS1}. For results on direct EM-scattering from screens see \cite{BC}, and for earlier results on inverse scattering see \cite{CCD}. For determination of a screen in dimensions three or higher from the full Dirichlet--Neumann map (including stability results), see \cite{AS} and references therein.

\section{Acknowledgements}

The authors are grateful to Kari Astala, who drew their attention to
Carathéodory's prime end theory. The work was supported by the
Estonian Research Council (grant number PRG832) and the Research
Council of Finland through the Flagship of Advanced Mathematics for
Sensing, Imaging and Modelling (decision number 359183).


\begin{thebibliography}{99}

\bibitem{AK} Akduman, I. and R. Kress: Electrostatic Imaging Via Conformal Mapping, {\em Inverse Problems}, Vol 18(6) (2002), 1659--1672, \url{http://dx.doi.org/10.1088/0266-5611/18/6/315}.

\bibitem{AS}  Alessandrini, G. and E. Sincich: Cracks with impedance; Stable
  determination from boundary data. {\em Indiana University Mathematics
  Journal}, 62(3), 947--989 (2013)

\bibitem{AH} C.~J.~S.~Alves and T.~Ha-Duong, {\em On inverse
  scattering by screens},  Inverse Problems , {\bf 13} (1997), no.~5,
  1161--1176.

\bibitem{BPS} Bl{\aa}sten, E., P\"aiv\"arinta, L. and S. Sadique: {\em Unique Determination of the Shape of a Scattering Screen from a Passive Measurement}, to appear.

\bibitem{BPS1} Bl{\aa}sten, E., P\"aiv\"arinta. L. and S. Sadique: The Fourier, Hilbert, and Mellin Transforms on a
Half-Line, {\em SIAM J. Math. Anal.} 55 (2023), no. 6, 7529--7548.

\bibitem{BC} Buffa, A. and S.H. Christiansen  {\em The electric field integral equation on Lipschitz
screens: definitions and numerical approximation}, Numer. Math, (2003) 94: 229--267

\bibitem{CCD}  Cakoni, F.,  Colton, D. and E. Darrigrand: The inverse electromagnetic scattering problem for screens
{\em Inverse Problems} 19 (2003), no. 3, 627--642.

\bibitem{CK-II} Colton, D. and R. Kress: {\em Inverse Acoustic and Electromagnetic Scattering Theory}, Springer (1992)

\bibitem{FV} Friedman, A. and M. Vogelius, M. : Determining Cracks By Boundary
  Measurements, {\em Indiana University Mathematics Journal}, 38(3), 527--556 (1989).

\bibitem{K-integral}  Kress, R.: {\em Linear Integral Equations}, Springer, (1989).

\bibitem{K-fixed-freq} Kress, R.: Fr\'echet differentiability of the far field operator for scattering from a crack, {\em Journal of Inverse and Ill-Posed Problems}, Vol 3(4) (1995), 305--313, \url{http://dx.doi.org/10.1515/jiip.1995.3.4.305}.

\bibitem{K-analytic} Kress, R.: Inverse scattering from an open arc, {\em Mathematical Methods in the Applied Sciences}, Vol 18(4) (1995), 267--293, \url{http://dx.doi.org/10.1002/mma.1670180403}.

\bibitem{LP} Lax, P. and R. Phillips: {\em Scattering Theory, Revised Edition}, Academic Press (1989)

\bibitem{McL} McLean, W.: {\em Strongly Elliptic Systems and Boundary Integral Equations}, Cambridge University Press (2000)

\bibitem {NPT} Nachman, A., P\"aiv\"arinta, L. and A. Teiril\"a: On Imaging Obstacles Inside Inhomogeneous Media. {\em J.Functional Analysis}, Vol 252(2) (2007), 490-516.

\bibitem{OPS} Ola, P., P\"aiv\"arinta. L. and S. Sadique: Unique Determination of a Planar Screen in Electromagnetic Inverse Scattering, {\em Mathematics}, 2023, 11(22), 4655

\bibitem{St} Stephan, E.: Boundary Integral Equations for Screen Problems in \(\R^3\), {\em Integrals Equations Operator Th} 10 (1987)

\bibitem{T2} Taylor, M.: {\em Partial Differential Equations II - Qualitative Studies of Linear Equations},  Springer (1996)

\end {thebibliography}
\end{document}